\documentclass[a4paper,thmsa]{article}
\usepackage{graphicx}
\usepackage{amsmath}
\usepackage{amsfonts}
\usepackage{amssymb}

\newtheorem{theorem}{Theorem}

\newtheorem{lemma}{Lemma}

\newenvironment{proof}[1][Proof]{\textbf{#1.} }{\ \rule{0.5em}{0.5em}}
\numberwithin{equation}{section}

\def\Vo{\mathrm{Var}_\omega}

\begin{document}

\title{Simple Transient Random Walks in One-dimensional Random
Environment: the Central Limit Theorem}
\author{Ilya Ya.
Goldsheid \\ School of Mathematical Sciences\\ Queen Mary,
University of London\\ London E1 4NS, Great Britain\\ email:
I.Goldsheid@qmul.ac.uk}
\date{May 30, 2006}
\maketitle
\begin{abstract}\noindent
We consider a simple random walk (dimension one, nearest neighbour
jumps) in a quenched random environment. The goal of this work is
to provide sufficient conditions, stated in terms of properties of
the environment, under which the Central Limit Theorem (CLT) holds
for the position of the walk. Verifying these conditions leads to
a complete solution of the problem in the case of independent
identically distributed environments as well as in the case of
uniformly ergodic (and thus also weakly mixing) environments.


\end{abstract}

\section{Introduction} The study of the asymptotic behaviour of
random walks (RW) in random environment (RWRE) has been started
more than thirty years ago. The first mathematical results were
obtained in the pioneering papers by M. Kozlov \cite{Koz}, Solomon
\cite{So}, and Kesten--M. Kozlov--Spitzer \cite{KKS}. The
asymptotic behaviour of a RW in annealed environments has been
described in \cite{KKS} in detail for all regimes except the
recurrent one. Sinai \cite{S} has completed this description by
discovering the $\log^2$ law in the recurrent case. Recently, the
results of \cite{KKS} were extended by
Mayer-Wolf--Roitershtein--Zeitouni to the case of Markovian
environments \cite{WRZ}.

But the question about the asymptotic behaviour of RW in a
quenched (frozen) environment remains largely open. However, it
has to be mentioned that Alili \cite{Al} proved the Central Limit
Theorem (CLT) for a random walk in a quasi-periodic environment
with very special additional properties.

The aim of this work is to prove that under certain sufficient
conditions (which are often also necessary) the Central Limit
Theorem (CLT) holds for simple (one-dimensional with nearest
neighbour jumps) random walks in a typical quenched environment.

Traditionally, the simple random walk is characterized by two
quantities: the hitting time $T(n)$ of site $n$  and the position
of the walk $X(t)$ at time $t$. One usually starts with the study
of the asymptotic behaviour of $T(n)$ as $n\to\infty$ and then
'translates' the results of this study into results for the
asymptotic behaviour of $X(t)$ as $t\to\infty$. Hitting times are
easy to control due to the fact (used already in \cite{So}) that,
in this model, they can be presented as sums of independent random
variables. In \cite{Al} the CLT for hitting times has been proved
for RW's in ergodic environments. The proof of this fact is given
below (Theorem \ref{ThMain1}) mainly because it is used in the
proof of Theorem \ref{ThMain2}. Our main results are concerned
with a less simple question about the position of the walk and are
as follows.

Theorem \ref{ThMain2}  reduces the question about the CLT for
$X(t)$ to a question about certain properties of the environment.
It also offers a choice of two random centerings for $X(t)$ which
are functions of the environment.

In Theorem \ref{ThMain3} we prove that independent identically
distributed (i.i.d.) random environments do have the properties
allowing to apply Theorem \ref{ThMain2} in this case. In fact it
can be shown (though we don't do it here) that environments
satisfying strong mixing conditions also have these properties.

Theorem \ref{ThMain4} provides a very short and simple proof of the
CLT for uniformly ergodic environments (in particular,
quasi-periodic environments). It thus addresses the other side of
the spectrum, as far as the mixing properties are concerned.

It should be emphasized that at present there is no proof of CLT
for a position of the walk which would work in a general ergodic
environment.

\smallskip The CLT in annealed setting is not discussed in this paper.
It can be derived from the quenched CLT in the case of
environments with strong mixing properties but it should be
stressed once gain that in the general ergodic setting even this
question remains opened both for the hitting times and the
position of the walk.

\smallskip Apart of the above there is the following reason for appearance of
this work. Our intention is to address the problem in the simplest
case since this is where the ideas can be best explained and the
proofs are short and transparent. The same approach, properly
adapted, works in a much more general case of a RWRE on a strip but
explaining it there is a much more technical matter.

\smallskip\noindent Since the problems considered in this work stem
directly from \cite{So, KKS}, we don't review the beautiful
development that followed the appearance of these papers.
Relatively recent and modern introductions to the subject as well
as comprehensive reviews can be found in \cite{Z}, \cite{BS1}, and
\cite{S5}.

\smallskip\noindent The article is organized as follows.

We start by describing the models considered in this work. We then
explain those results from \cite{So}, \cite{KKS} and \cite{Al}
which are relevant to this work. This is followed by statement of
our main results which are then proved in the next section.
Appendix contains several technical results some of which may be
new and some are unlikely to be new but are included mainly for
the sake of completeness.

I am grateful to E. Bolthausen and O. Zeitouni for valuable
discussions and suggestions.

This paper was about to be submitted when I learned that O.
Zeitouni \& J. Peterson obtained a result which is similar to the
one stated in Theorem \ref{ThMain3}.

\subsection{Description of the model.\label{s1.1}}

Let $(\Omega,\mathcal{F},\mathbb{P} ,\mathcal{T})$ be a dynamical
system with $\Omega=\{\omega\}$ denoting a set of elementary
events, $\mathcal{F}$ being a $\sigma$-algebra of subsets of
$\Omega$, $\mathbb{P}$ denoting a probability measure on
$(\Omega,\mathcal{F})$, and $\mathcal{T}: \Omega\mapsto\Omega$
being an invertible transformation of $\Omega$ preserving measure
$\mathbb{P}$. Next, let $p\,:\,\Omega\mapsto (0,1)$ be a
measurable real valued function on $\Omega$ such that
$0<p(\omega)<1$ for all $\omega\in\Omega$.

Put $p_n\equiv p_{n}(\omega)=p(\mathcal{T}^{n}\omega)$,
$q_n=1-p_n$, $-\infty<n<\infty$. For any such sequence $p_n$ we
shall now define a random walk $X(t;\omega,z)$ with discrete
integer valued time $t\ge0$. The phase space of the walk is a
one-dimensional lattice $\mathbb{Z}$ and $p_n,\ q_n$ are its
transition probabilities:
\begin{equation}
\mathcal{Q}_{\omega}(z_1,z_{2})\overset{\mathrm{def}}{=}\left\{
\begin{array}
[c]{ll}%
p_n & \mathrm{if\quad}z_1=n,\ z_{2}=n+1,\\
q_{n} & \mathrm{if\quad}z_1=n,\ z_{2}=n-1,\\
0 & \mathrm{otherwise.}%
\end{array}
\right. \label{1}%
\end{equation}
For any starting point $z\in\mathbb{Z}$ the probability law on the
space of trajectories is denoted by $Pr_{\omega,z}$ and is defined
    by its finite-dimensional distributions
\begin{equation}
Pr_{\omega,z}\left(  X(1)=z_{1},\ldots,X(t)=z_{t}\right) \overset
{\mathrm{def}}{=}\mathcal{Q}_{\omega}(z,z_{1})\mathcal{Q}_{\omega}(z_{1}%
,z_{2})\cdots\mathcal{Q}_{\omega}(z_{t-1},z_{t}).\label{StripRWRE}%
\end{equation}

We say that the sequence $p_n$  (or, equivalently, the $\omega$) is
the \textit{environment} or the \textit{random environment} of the
walk. The \textit{annealed probability measure} on the product of
the space  of environments $\Omega$ and the space of trajectories
$X(\cdot;\omega,z)$ starting form $z$ is a semi-direct product of
$\mathbb{P}\times Pr_{\omega,z}$, defined by
$\mathbb{P}(d\omega)Pr_{\omega ,z}(dX)$. We write $Pr_{\omega}$ for
$Pr_{\omega,z}$ and $X(\cdot)$ for $X(\cdot;\omega,z)$ when there is
no danger of confusion. It is useful to remember that, unless
explicitly stated otherwise, we always suppose that the environment
is quenched (frozen).

\smallskip The just described general class of models provides a natural
setting for Theorems \ref{ThMain1} and \ref{ThMain2}.

It has already been mentioned above that we shall consider two
sub-classes of this model. The so called {\it the i.i.d.
environments} form one of these sub-classes and arise when $p_n$
is a sequence of independent identically distributed random
variables.

A sub-class of random environments which we call {\it uniformly
ergodic environments} is obtained when the dynamical system has
very good ergodic properties which are usually combined with very
week mixing properties. It is convenient to give the precise
definition later but it is natural to mention here that a
quasi-periodic environment is also a uniformly ergodic
environment.

\subsection{Notations and assumptions.\label{s1.2}}

\smallskip\noindent
{\it 1. Hitting times.} Let $T_k(n)$ be the hitting time of site $n$
by a random walk $X(\cdot;\omega,k)$ starting from $k$:
 \[
 T_k(n)\overset{\mathrm{def}}{=}\inf\{ t:\ X(t;\omega,k)=n\}.
 \]
The notation $T(n)$ is reserved for the case $k=0$. We put
$\tau_k=T_k(k+1)$. The random variables $\tau_k$ are independent
when $\omega$ is fixed (with their distributions depending on $k$
and $\omega$). As in \cite{So}, we shall make use of the following
simple relation
\begin{equation}\label{To}
T_k(n)=\sum_{j=k}^{n-1}\tau_k.
\end{equation}
\smallskip\noindent
{\it 2. Expectations.} Throughout the paper $\mathbb{E}$  denotes
the expectation with respect to the measure $\mathbb{P}$. By
$\mathcal{E}_{\omega,z}$ we denote the expectation with respect to
the measure $Pr_{\omega,z}$; in those cases when the starting
point of the walk is clearly defined by the context we may use
$\mathcal{E}_{\omega}$ for $\mathcal{E}_{\omega,z}$ (e. g.
$\mathcal{E}_{\omega}\tau_k\equiv\mathcal{E}_{\omega,k-1}\tau_k
$). The notation $\Vo$ will be used for the variance of a random
variable calculated with respect to the measure $Pr_{\omega,z}$,
e. g. $\Vo(\tau_k)=\mathcal{E}_{\omega}(\tau_k-
\mathcal{E}_{\omega}\tau_k)^2$.

\smallskip\noindent
{\it 3. Main assumptions.} The following set of assumptions is
called Condition $\mathbf{C}$ and is supposed to be satisfied
throughout the paper:

\noindent \textbf{Condition $\mathbf{C}$}
\begin{description}
\item [C1]The dynamical system
$(\Omega,\mathcal{F},\mathbb{P},\mathcal{T}) $ is ergodic

\item[C2] \ \ \ \ \ $\mathbb{E}\log p_k^{-1}<\infty,\ \ \mathbb{E}%
\log(1-p_k)^{-1}<\infty $
\end{description}
A set of stronger assumptions called Condition $\mathbf{C'}$
consists of $\mathbf{C1}$, $\mathbf{C3}$ and $\mathbf{C4}$:
\begin{description}
\item[$\mathbf{C3}$] There is a $\gamma>2$ such that
\[
\mathbb{E}p_k^{-\gamma}<\infty,\ \ \mathbb{E}%
(1-p_k)^{-\gamma}<\infty
\]
\item[$\mathbf{C4}$]
$\limsup_{n\to\infty}(\mathbb{E}\prod_{j=1}^n(\frac{q_j}{p_j})^\gamma)^{\frac{1}{n}}<\infty$
\end{description}
{\it Remark.} Obviously, $\mathbf{C4}$ follows from $\mathbf{C3}$
if the environment is i.i.d. This is not true in general ergodic
setting.

\subsection{Preliminaries: transience, recurrence, linear growth.}
We say that a random walk (in a fixed environment $\omega$) is
transient to the right (transient to the left) if
\[
\lim_{t\rightarrow\infty}X(t)=\infty\ \ (\hbox{correspondingly}\ \
\lim_{t\rightarrow\infty}X(t)=-\infty\ ).
\]
We shall now quote several statements from \cite{So} in a form
which suits us best. Let us put
\begin{equation}
A_j\overset{\mathrm{def}}{=}\frac{q_j}{p_j}, \ \hbox{ and }\
\lambda\overset{\mathrm{def}}{=} \mathbb{E}\ln A_j
\label{A}%
\end{equation}
(it is clear that $\lambda$ does not depend on $j$).

\smallskip\noindent {\it The recurrence and transience criteria} for
our random walk are given by the following result from \cite{So}.
\begin{theorem}
\label{ThLyap} Suppose that Condition $\mathbf{C}$ is satisfied.
Then

\noindent (i) $\lambda<0$ implies for $\mathbb{P}$-a.e. environment
$\omega$ that $X_\omega$ is transient to the right. Symmetrically,
$\lambda>0$, implies for $\mathbb{P}$-a.e. environment $\omega$ that
$X_\omega$ is transient to the left.

\noindent (ii) $\lambda=0$ if and only if $X_\omega$ is recurrent
for $\mathbb{P}$-a.e. $\omega$, that is
\[
\limsup_{t\rightarrow\infty}X(t)=+\infty,\ \ \ \liminf_{t\rightarrow\infty}%
X(t)=-\infty \hbox{ $Pr_{\omega}$-almost surely}.\]
\end{theorem}
\noindent\textbf{From now on} we consider only those RWRE which
are transient to the right, that is $\lambda <0$. To state further
results it is convenient to define a function $r(\kappa)$
depending on a parameter $\kappa\in [0,\,\gamma]$, where $\gamma$
is the same as in $\mathbf{C2'}$, namely:
\begin{equation}
r(\kappa)\overset{\mathrm{def}}{=}
\limsup_{n\to\infty}\left(\mathbb{E}\prod_{j=1}^nA_j^\kappa\right)^{\frac{1}{n}}.
\label{r}%
\end{equation}
This function is a simple generalization of the one first
considered in \cite{KKS} (see also \cite{WRZ} where $\log
r(\kappa)$ has been studied). If the $p_n$'s are i.i.d. random
variables then of course
\[
r(\kappa)= \mathbb{E}\left(A_0 \right)^\kappa.
\]
As has been shown in \cite{KKS} and in \cite{WRZ}, the asymptotic
behaviour of the RWRE can be characterized in terms properties of
$r(\kappa)$ which are well worth of being studied. However, for
the purposes of this work, we only need the following simple

\begin{lemma}
\label{simple} Suppose that Condition $\mathbf{C}'$ is satisfied.
Then the function $\ln r(\kappa)$ is continuous and convex on
$[0,\, \gamma)$.
\end{lemma}
The proof of this Lemma is given in the Appendix.

Let us define a function which plays a very important role in this
paper (as it did already in \cite{So}): for a fixed environment
$\omega$ put
\begin{equation}\label{mu0}
\mu_0(\omega)\overset{\mathrm{def}}{=}1+2A_0+\ldots+2A_0A_{-1}\ldots
A_{-j}+\ldots\equiv 1+2\sum_{j=0}^\infty \prod_{i=-j}^0A_i,
\end{equation}
and
\begin{equation}\label{mu}
\mu_k(\omega)\overset{\mathrm{def}}{=}\mu_0(\mathcal{T}^k\omega)
\equiv 1+2\sum_{j=0}^\infty \prod_{i=k-j}^kA_i.
\end{equation}
The probabilistic meaning of $\mu_k(\omega)$ is explainne by the
following
\begin{lemma}(\cite{So}, \cite{Z})\label{Etau} If $\lambda<0$ then $\mu_k$
is finite for $\mathbb{P}$-almost all $\omega$ and
$\mathcal{E}_\omega\tau_k=\mu_k$.
\end{lemma}
Let us note first that $\mu_k(\omega)$ has the following property:
\begin{equation}\label{Emu}
\hbox{if $r(\kappa)<1$ then there is a $\delta>0$ such that
$\mathbb{E}\mu_k^{\kappa+\delta}<\infty $}.
\end{equation}
It is easy to see that (\ref{Emu}) holds for any $\kappa>0$ but
since we need it when $\kappa\ge1$, it shall be explained only in
this case (the other one is even simpler).

Namely, according to Lemma \ref{simple}, $r(\cdot)$ is a
continuous function. We thus can choose $\delta>0$ and such that
$r(\kappa+\delta)<1$. Consider the Banach space
$L_{\kappa+\delta}(\Omega)$ with $||Y||\overset{\mathrm{def}}{=}
(\mathbb{E}|Y|^{\kappa+\delta})^\frac{1}{\kappa+\delta}$ for any
function $Y\in L_{\kappa+\delta}(\Omega)$. We then have:
\[
||\mu_0||_{\kappa+\delta}\le 1+2\sum_{j=0}^\infty ||\prod_{i=
-j}^0A_i||_{\kappa+\delta}<\infty,
\]
and this proves (\ref{Emu}) (remember that $||\mu_0||=||\mu_k||$).

In particular if $r(1)<1$ then
\begin{equation}\label{Emu1}
\mu\overset{\mathrm{def}}{=}\mathbb{E}\mu_k \le
(\mathbb{E}\mu_k^{1+\delta})^\frac{1}{1+\delta}<\infty\ \  \hbox{
for some \ $\delta>0$}.
\end{equation}

\smallskip\noindent {\it The quenched Law of Large Numbers} has
been proved in \cite{So} for i.i.d. environments and the same
proof works in general ergodic setting (see \cite{Al} or \cite{Z}
for more detailed explanations).

\begin{theorem} \label{ThMain} Suppose that Condition
$\mathbf{C}'$ is satisfied and that $\lambda<0$. Then:

\noindent (i) $r(1)<1$  implies that for $\mathbb{P}$-a.e.
environment $\omega$ with $Pr_{\omega}$-probability 1
\begin{equation}\label{1.9}
\lim_{n\rightarrow\infty}\frac{T(n)}{n}=\mu<\infty\ \hbox{ and }\
\lim_{t\rightarrow\infty}\frac{X(t)}{t}=\mu^{-1}>0,
\end{equation}
\noindent (ii) $r(1)>1$ implies that for $\mathbb{P}$-a.e.
environment $\omega$ with $Pr_{\omega}$-probability 1
\begin{equation}\label{1.10}
\lim_{n\rightarrow\infty}\frac{T(n)}{n}=\infty \ \hbox{ and }\
\lim_{t\rightarrow\infty}\frac{X(t)}{t}=0.
\end{equation}
\noindent (iii) If the environment is i.i.d. then (\ref{1.10})
holds also for $r(1)=1$.
\end{theorem}
{\it Remark.} If the environment is i.i.d. then a straightforward
calculation leads to an explicit formula for $\mu$ (known since
\cite{So}): $ \mu=\frac{1+r(1)}{1-r(1)}$.

We finish this section by defining uniformly ergodic environments.

\noindent\textbf{Definition 1} \textit{ Let
$f:\Omega\mapsto\mathbb{R}$ be an $\mathcal{F}$-measurable
function on $\Omega$. We say that the transformation $\mathcal{T}$
is $f$-uniformly ergodic if
\begin{equation}\label{uniferg}
\left|n^{-1}\sum_{j=1}^{n}f(\mathcal{T}^j\omega)-\mathbb{E}f\right|\le\varepsilon_n
\end{equation}
where the sequence $\varepsilon_n$ does not depend on $\omega$ and
$\lim_{n\to\infty}\varepsilon_n=0$.}

\noindent \textit{We say that a random environment is
\textbf{uniformly ergodic} if $\mathcal{T}$ is $\mu_0$-uniformly
ergodic. }

One of the simplest uniformly ergodic environments is generated by
a quasi-periodic dynamical system with $\Omega=[0,1]$,
$\mathcal{T}(\omega)=(\omega +\alpha)(\mod1)$, where
$\alpha\in[0,1]$ is an irrational number. If the function
$p(\cdot)$ is continuous on $[0,1]$, $p(0)=p(1)$, and that $
\lambda=\int_0^1\ln\frac{1-p(\omega)}{p(\omega)}d\omega <0 $ then
also $\mu_0(\omega)$ a continuous function on $[0,1]$ and the
uniform ergodicity of this environment follows.

Let us explain this statement in a more general setting. Suppose
that $\mathcal{T}$ is a continuous homeomorphism of a compact
metric space $\Omega$ and that $\mathbb{P}$ is its unique
invariant measure. Suppose also that the function $p(\cdot)$ is
continuous. It is then easy to see that
$r(\kappa)=e^{\kappa\lambda}$. Indeed,
\[
(\mathbb{E}\prod_{j=1}^nA_j^\kappa)^{\frac{1}{n}}=
(\mathbb{E}e^{\kappa\sum_{j=1}^n\ln A_j})^{\frac{1}{n}}=
(\mathbb{E}e^{\kappa n(\lambda+\epsilon_n)})^{\frac{1}{n}},
\]
where $|\epsilon_n(\omega)|\le \varepsilon_n$. Hence $e^{\kappa
(\lambda-\varepsilon_n)}\le r(\kappa)\le e^{\kappa
(\lambda+\varepsilon_n)}$ and the statement follows. If now
$\lambda<0$, then series (\ref{mu0}) converges uniformly in
$\omega\in\Omega$ and hence $\mu_0$ is a continuous function on
$\Omega$. The latter in  turn implies uniform ergodicity of the
environment.

\section{Main Results}
In order to state the central limit theorem for $T(n)$ one has to
know the variance of this random variable. It turns out that in
the case of the simple walk an explicit expression for the
variance can be found and the calculations are not complicated.
Formula (\ref{2.12}) has been obtained in \cite{Al} where
branching processes are used for its derivation. We use a
different approach which works also for more general models
(\cite{G}).
\begin{lemma} \label{var} Suppose that $\lambda<0$. Then for
$\mathbb{P}$-almost every $\omega$ the variance of $T(n)\equiv
T(n;\omega)$ is finite and is given by
\begin{equation}\label{2.11}
\Vo(T(n)) = \sum_{k=0}^{n-1}\sigma_k^2(\omega),
\end{equation}
where
\begin{equation}\label{2.12}
\sigma_k^2(\omega)\overset{\mathrm{def}}{=}
\Vo(\tau_k)=\sum_{j=0}^\infty
p_{k-j}^{-1}(\mu_{k-j-1}+1)^2\prod_{i=k-j}^kA_i.
\end{equation}
If in addition $r(2)<1$, then
\begin{equation}\label{2.13}
\sigma^2\equiv\mathbb{E}\Vo(\tau_k)<\infty.
\end{equation}
If $r(2)<1$ and the environment is i.i.d. then
\begin{equation}\label{2.14}
\sigma^2=\frac{4(r(1)+r(2))(1+r(1)^2)}{(1-r(1))^2(1-r(2))}.
\end{equation}
\end{lemma}
The proof of Lemma \ref{var} is given in appendix.

\smallskip\noindent Denote
\begin{equation}\label{H}
H(n,\omega)=\sum_{k=0}^{n-1}\mu_k\equiv\mathcal{E}_\omega(T(n)),
\end{equation}
where the last equality follows from (\ref{To}) and Lemma
\ref{Etau}. We often write $H(n)$  for $H(n,\omega)$. It is clear
from \ref{H} that $H(n)$ is the natural centering in the CLT for
$T(n)$. It turns out that centerings for $X(t)$ can too be
expressed, with a varying degree of explicitness, in terms of the
function $H(\cdot)$.

In the sequel we denote $\lfloor
y\rfloor\overset{\mathrm{def}}{=}$ integer part of $y$, where $y$
is any real number. We also use the following convention about
summations. For any real numbers $b_1,\ b_2$ and a sequence $d_k,\
-\infty<k<\infty$,
\begin{equation}\label{notations2}
\sum_{k=b_1}^{b_2}d_k\overset{\mathrm{def}}{=} \sum_{k=\lfloor
b_1\rfloor}^{\lfloor b_2\rfloor}d_k \overset{\mathrm{def}}{=}
-\sum_{k=\lfloor b_2\rfloor}^{\lfloor b_1\rfloor}d_k.
\end{equation}
In particular for $y\ge0$ we put $
H(y)\overset{\mathrm{def}}{=}H(\lfloor y\rfloor)$.

\smallskip\noindent\textbf{Definition 2.} \textit{The function
\begin{equation}\label{b}
b(t;\omega)\overset{\mathrm{def}}{=}2\mu^{-1}t -\mu^{-1}
H(\mu^{-1}t,\omega)
\end{equation}
is said to be the explicit centering for $X(t)$. The integer
valued function $\tilde b(t;\omega)$ such that
\begin{equation}\label{b1}
H(\tilde{b}(t;\omega)) \equiv
\sum_{k=0}^{\tilde{b}(t;\omega)-1}\mu_k\le t<
\sum_{k=0}^{\tilde{b}(t;\omega)}\mu_k\equiv
H(\tilde{b}(t;\omega)+1).
\end{equation}
is said to be the implicit centering for $X(t)$.}

\smallskip\noindent It is easy to see that
\begin{equation}\label{b2}
b(t;\omega)=\mu^{-1}t -\mu^{-1}\sum_{k=0}^{\mu^{-1}t-1}
(\mu_k-\mu) + O_1(1),
\end{equation}
where $O_1(1)=\mu^{-1}t-\lfloor\mu^{-1}t\rfloor<1$.

\smallskip\noindent\textbf{For the rest} of the paper we suppose that
\[
\hbox{ Condition $\mathbf{C}'$ is satisfied and $r(2)<1$.}
\]
Put $\Phi(x)=\frac{1}{\sqrt{2\pi}}\int_{-\infty}^x
e^{-\frac{u^2}{2}}du$. We shall prove the following statements.
\begin{theorem}
\label{ThMain1} (\cite{Al}) For for $\mathbb{P}$-almost every
environment $\omega$
\begin{equation}\label{clt1}
Pr_\omega\left\{\frac{T(n)- H(n)}{\sqrt{n}\sigma}<x\right\} \to
\Phi(x)\ \ \hbox{ uniformly in $x$ as $n\to\infty$}.
\end{equation}
\end{theorem}
{\it Remark.} Uniform convergence in (\ref{clt1}) is used in the
proof of Theorem \ref{ThMain2}.
\begin{theorem}
 \label{ThMain2} Suppose that at least one
of the following two relations holds:
\begin{equation}\label{rel1}
\lim_{t\to\infty}t^{-\frac{1}{2}}
\sum_{k=\mu^{-1}t}^{b(t)+\sqrt{t}\sigma^*x}(\mu_k-\mu)=0 \ \
\hbox{ with $\mathbb{P}$-probability 1 for any real $x$ }
\end{equation}
\begin{equation}\label{rel2}
\lim_{t\to\infty}t^{-\frac{1}{2}}
\sum_{k=\tilde{b}(t)}^{\tilde{b}(t)+\sqrt{t}\sigma^*x}(\mu_k-\mu)=0
\ \ \hbox{ with $\mathbb{P}$-probability 1 for any real $x$}.
\end{equation}
Then for $\mathbb{P}$-almost every environment $\omega$
\begin{equation}\label{clt2}
\lim_{t\rightarrow\infty}Pr_\omega\left\{\frac{X(t)-
\bar{b}(t)}{\sqrt{t}\sigma^*} \le x\right\}= \Phi(x),
\end{equation}
where ${\sigma^*}^2=\mu^{-3}\sigma^2$, the convergence in
\ref{clt2} is uniform in $x$ and
\[
\bar{b}(t)=
\begin{cases}
{b}(t) & \hbox{ if (\ref{rel1}) holds}\\
\tilde{b}(t) & \hbox{ if (\ref{rel2}) holds}\\
\end{cases}
\]
($\bar{b}(t)$ can be equal to any of the two if both (\ref{rel1})
and (\ref{rel2}) hold).
\end{theorem}
{\it Remarks.} 1. Note that in (\ref{rel1}) and in (\ref{rel2})
the summation is carried out within random limits.

\noindent 2. The explicit form in which $b(t)$ is given by
(\ref{b}) allows one to state condition (\ref{rel1}) in the
following equivalent form:
\[
\lim_{t\to\infty}t^{-\frac{1}{2}}
\sum_{k=t}^{t-\mu^{-1}\sum_{j=o}^{t-1}(\mu_j-\mu)+\sqrt{t}y}(\mu_k-\mu)=0
\ \ \hbox{ with $\mathbb{P}$-probability 1 for any real $y$}.
\]
However, (\ref{rel1}) is in fact a good approximation for
(\ref{rel2}) in the case of environments with sufficiently strong
mixing properties. Besides, it is also more convenient to use it
in the proof of Theorem \ref{ThMain2}.

\smallskip\noindent We finish this section by stating two theorems which
demonstrate the usefulness of conditions (\ref{rel1}) and
(\ref{rel2}).
\begin{theorem}
\label{ThMain3} In the i.i.d.random environment (\ref{rel1}) holds
and thus for $\mathbb{P}$-almost every environment $\omega$
\begin{equation}\label{clt3}
\lim_{t\rightarrow\infty}Pr_\omega\left\{\frac{X(t)-
b(t)}{\sqrt{t}\sigma^*} \le x\right\}= \Phi(x)
\end{equation}
with convergence in (\ref{clt3}) being uniform in $x$.
\end{theorem}
\begin{theorem}
\label{ThMain4} Suppose that the environment is uniformly ergodic.
Then (\ref{rel2}) holds and thus for $\mathbb{P}$-almost every
environment $\omega$
\begin{equation}\label{clt4}
\lim_{t\rightarrow\infty}Pr_\omega\left\{\frac{X(t)-
\tilde{b}(t)}{\sqrt{t}\sigma^*} \le x\right\}= \Phi(x)
\end{equation}
and convergence in (\ref{clt4}) is uniform in $x$.
\end{theorem}

\section{Proofs}
{\bf Proof of Theorem \ref{ThMain1}.} Proving (\ref{clt1})
essentially means proving a CLT for the sum of independent random
variables $\tau_k$. Indeed, since
\begin{equation}\label{3.0}
\frac{T(n)- H(n,\omega)}{\sqrt{n}\sigma}=\frac{T(n)-
\mathcal{E}_\omega(T(n))}{\sqrt{\Vo(T(n))}}\sqrt{\frac{\Vo(T(n))}{n\sigma^2}},
\end{equation}
it is enough to check that for $\mathbb{P}$-almost all $\omega$
\begin{equation}\label{3.1}
\frac{\Vo(T(n))}{n\sigma^2}\to1
\end{equation}
and that CLT holds for $T(n)$.  Relation (\ref{3.1})  follows from
(\ref{2.11}), (\ref{2.13}), and the Birkhoff ergodic theorem.
Next, for those $\omega$ for which (\ref{3.1}) holds, also
\[
\max_{0\le k\le n-1}\sigma_k^2/\Vo(T(n))\to 0 \hbox{ as $n\to
\infty$}.
\]
This in turn is well known to imply that the Lindeberg's
conditions for the CLT for sums of non-identically distributed
random variables $\tau_k$ holds. Theorem \ref{ThMain1} is proved. $\Box$

\smallskip\noindent{\bf Proof of Theorem \ref{ThMain2}.}
As usual with CLT's, it is sufficient to prove (\ref{clt2}) for
every fixed value of $x$. Thus, for the duration of the proof $x$,
is considered to be a fixed parameter. The proof of (\ref{clt2})
will be split into three parts.

\noindent{\it Part 1: approximating $X(t)$ by $n_t$.} For any
integer time $t$ let $n_t$ be a positive random integer such that
\begin{equation}\label{p1.1}
T(n_t)\le t < T(n_t+1).
\end{equation}
{\it Remark.} This definition of $n_t$ has been used already in
\cite{So} in the proof of the Law of Large Numbers cited above.

Since $|X(t)-X(t')|\le|t-t'| $ and since $X(T(n_t))=n_t$ and
$X(T(n_t+1))=n_t+1$ (by the definition of $T(n)$), it follows that
\begin{equation}\label{p1.2}
|X(t)-n_t|=|X(t)-X(T(n_t))|\le t -
T(n_t)<T(n_t+1)-T(n_t)=\tau_{n_t}.
\end{equation}
Hence
\begin{equation}\label{p1.3}
t^{-\frac{1}{2}}|X(t)-n_t|<t^{-\frac{1}{2}}\tau_{n_t}.
\end{equation}
The sequence $\tau_k$ forms a stationary process in annealed
environment and since  $\mathbb{E}\left(\mathcal{E}_\omega
\tau_k^2\right)<\infty$ we have that
$t^{-\frac{1}{2}}\tau_{n_t}\to0$ as $t\to\infty$ with
$\mathbb{P}\times Pr_\omega$-probability 1 which in turn implies
that it holds for $\mathbb{P}$ -almost every $\omega$ with
$Pr_\omega$-probability 1. This implies that proving (\ref{clt2}) is
equivalent to proving that
\begin{equation}\label{p1.4}
\lim_{t\rightarrow\infty}Pr_\omega\left\{\frac{n_t-
\bar{b}(t)}{\sqrt{t}\sigma^*} \le x\right\}= \Phi(x).
\end{equation}
{\it Part 2: proof for the case when (\ref{rel1}) holds.} We need
a simple (but very useful) identity. Namely, it follows from
(\ref{p1.1}) and monotonicity of the function $T(\cdot)$ that for
any $y\ge 0$ the following two events coincide:
\begin{equation}\label{p1.5}
\{n_t\le y\}=\{T(y+1)> t\},
\end{equation}
where as before $T(y+1)\equiv T(\lfloor y\rfloor+1)$. This
identity is a slight modification of the one  which has been often
used in the context of RWRE at least since the appearance of paper
\cite{KKS}.

Hence, for sufficiently large values of $t$ we can write
\begin{equation}\label{p1.6}
\begin{aligned} &
Pr_\omega\left\{\frac{n_t-b(t)}{\sqrt{t}\sigma^*} \le x\right\}=\\
&Pr_\omega\left\{T(b(t)+\sqrt{t}\sigma^*x+1) > t\right\} \equiv
Pr_\omega\left\{T(B(t)+1) > t\right\},
\end{aligned}
\end{equation}
where $ B(t)\overset{\mathrm{def}}{=}b(t)+\sqrt{t}\sigma^*x.$ It
is natural to use the fact that, for a typical fixed $\omega$,
$T(B(t)+1)$ is an asymptotically normal random variable. Let us
take a closer look at the parameters of this random variable.
First of all it follows from (\ref{b}) and the Birkhoff ergodic
theorem that
\begin{equation}\label{U1}
\lim_{t\to\infty}t^{-1}B(t)=\mu^{-1}\ \hbox{ for
$\mathbb{P}$-almost every }\ \omega.
\end{equation}
Relation (\ref{U1}) and  the Birkhoff ergodic theorem imply that
for $\mathbb{P}$-almost all $\omega$
\begin{equation}\label{U2}
t^{-1}\Vo(T(B(t)+1))= t^{-1}\sum_{k=0}^{B(t)}\sigma_k^2(\omega)
\to\mu^{-1}\sigma^2\ \ \hbox{ as }\ \  t\to\infty.
\end{equation}
Finally,
\begin{equation}\label{p1.7}
\begin{aligned} &
\mathcal{E}_\omega \left(B(t)+1)\right) =\sum_{k=0}^{B(t)}\mu_k=
B(t)\mu+ \sum_{k=0}^{B(t)}(\mu_k-\mu)+O_2(1),
\end{aligned}
\end{equation}
where $O_2(1)=\mu +\mu(\lfloor B(t) \rfloor -B(t))\le\mu$.
Returning to the original expression for  $B(t)$ and
simultaneously replacing $b(t)$ in the right hand side of
(\ref{p1.7}) by its expression from (\ref{b2}) leads to
\begin{equation}\label{p1.8}
\begin{aligned}
& \mathcal{E}_\omega \left(T(B(t)+1)\right)
=\\
& t+\sqrt{t}\sigma^* x \mu -\sum_{k=0}^{\mu^{-1}t}(\mu_k-\mu)
+\sum_{k=0}^{b(t)+\sqrt{t}\sigma^*x}(\mu_k-\mu)+O_3(1)=\\
& t+\sqrt{t}\sigma^* x \mu
+\sum_{k=\mu^{-1}t}^{b(t)+\sqrt{t}\sigma^*x}(\mu_k-\mu) +O_3(1),
\end{aligned}
\end{equation}
where
$O_3(1)=\mu(\mu^{-1}t-\lfloor\mu^{-1}t\rfloor)+O_2(1)\le2\mu$.
Putting
\[
\mathcal{G}(t)=\frac{T(B(t)+1)-\mathcal{E}_\omega
\left(T(B(t)+1)\right)}{\sqrt{t}\sigma^* \mu}
\]
we can present the right hand side of (\ref{p1.6}) as
\begin{equation}\label{p1.9}
\begin{aligned}
& Pr_\omega\left\{T(b(t)+\sqrt{t}\sigma^*x) > t\right\}=\\
& Pr_\omega\left\{ \mathcal{G}(t)> -x
-t^{-\frac{1}{2}}(\sigma^*\mu)^{-1}
\sum_{k=\mu^{-1}t}^{b(t)+\sqrt{t}\sigma^*x}(\mu_k-\mu) +
o(t^{-\frac{1}{2}}) \right\}.
\end{aligned}
\end{equation}
But, according to (\ref{rel1}), we have for
$\mathbb{P}$-almost every $\omega$:
\begin{equation}\label{p1.10}
\lim_{t\to\infty}t^{-\frac{1}{2}}
\sum_{k=\mu^{-1}t}^{b(t)+\sqrt{t}\sigma^*x}(\mu_k-\mu)=0,
\end{equation}
and also, because of (\ref{clt1}), we have that for
$\mathbb{P}$-almost every $\omega$ the sequence $\mathcal{G}(t)$
converges in distribution to a standard normal random variable.
Hence
\[
\lim_{t\to\infty} Pr_\omega\left\{ \mathcal{G}(t)> -x
-t^{-\frac{1}{2}}(\sigma^*\mu)^{-1}
\sum_{k=\mu^{-1}t}^{b(t)+\sqrt{t}\sigma^*x}(\mu_k-\mu)
+o(t^{-\frac{1}{2}}) \right\}=\Phi(x)
\]
which proves (\ref{clt2}).

\smallskip\noindent{\it Part 3: proof in the case when (\ref{rel2}) holds.}
The proof goes along the same lines as in Part 2 with the natural
replacement of $b(t)$ by $\tilde{b}(t)$. On the other hand, subtle
differences appear at the end of the proof; this is why it may be
useful to give a brief outline of it here.

As in Part 2, it follows from (\ref{p1.5}) that
\begin{equation}\label{p2.6}
Pr_\omega\left\{\frac{n_t-\tilde{b}(t)}{\sqrt{t}\sigma^*} \le
x\right\}= Pr_\omega\left\{T(\tilde{B}(t)+1) > t\right\},
\end{equation}
where $\tilde{B}(t)\overset{\mathrm{def}}{=}
\tilde{b}(t)+\sqrt{t}\sigma^*x.$ Also for $\mathbb{P}$-almost all
$\omega$
\begin{equation}\label{pU1}
\lim_{t\to\infty}t^{-1}B(t)=\mu^{-1}\ \hbox{ and }\ \
\lim_{t\to\infty}t^{-1}\Vo(T(B(t)+1))= \mu^{-1}\sigma^2.
\end{equation}
Next
\begin{equation}\label{p2.8}
\begin{aligned}
& \mathcal{E}_\omega \left(T(\tilde{B}(t)+1)\right)
=\sum_{k=0}^{\tilde{b}(t)+\sqrt{t}\sigma^*x}\mu_k=
\sum_{k=0}^{\tilde{b}(t)-1}\mu_k+
\sum_{k=\tilde{b}(t)}^{\tilde{b}(t)+\sqrt{t}\sigma^*x}\mu_k\\
& t+\sqrt{t}\sigma^* x \mu
+\sum_{k=\tilde{b}(t)}^{\tilde{b}(t)+\sqrt{t}\sigma^*x}(\mu_k-\mu)
+O_4(1),
\end{aligned}
\end{equation}
where the second line in (\ref{p2.8}) follows from the definition
of $\tilde{b}(t)$ (see (\ref{b1})) and
\[
O_4(1)=t- \sum_{k=0}^{\tilde{b}(t)-1}\mu_k +\mu(1+
\sqrt{t}\sigma^*x-
\lfloor\sqrt{t}\sigma^*x\rfloor)\le\mu_{\tilde{b}(t)}+2\mu.
\]
Let us denote $ \tilde{\mathcal{G}}(t)=(\sqrt{t}\sigma^*
\mu)^{-1}\left(T(\tilde{B}(t)+1)-\mathcal{E}_\omega
\left(T(\tilde{B}(t)+1)\right)\right) $. We can then right that
\begin{equation}\label{p1.11}
\begin{aligned}
& Pr_\omega\left\{T(\tilde{b}(t)+\sqrt{t}\sigma^*x) \ge t\right\}=\\
& Pr_\omega\left\{ \tilde{\mathcal{G}}(t)\ge -x
-t^{-\frac{1}{2}}(\sigma^*\mu)^{-1}
\sum_{k=\tilde{b}(t)}^{\tilde{b}(t)+\sqrt{t}\sigma^*x}(\mu_k-\mu)
+t^{-\frac{1}{2}}O_5(1))\right\},
\end{aligned}
\end{equation}
Where $O_5(1)$ is proportional to $O_4(1)$. We note that $
t^{-\frac{1}{2}}\mu_{\tilde{b}(t)}\to 0 $ with
$\mathbb{P}$-probability 1 because $\mu_k^2$ is a stationary
sequence with $\mathbb{E}\mu_k^2<\infty$ (see (\ref{Emu})). This
together with (\ref{rel2}) and the asymptotic normality of
$\tilde{\mathcal{G}}(t)$ finishes the proof. $\Box$

\smallskip\noindent{\bf Proof of Theorem \ref{ThMain3}.} According to
Theorem \ref{ThMain2} we only have to check that for i.i.d.
environments (\ref{rel1}) holds true. In fact, we shall prove that
i.i.d. environments satisfy (\ref{rel3}) which is slightly stronger
than (\ref{rel1}). To explain the last statement let us put
\begin{equation}\label{H1}
\mathcal{H}(n,\omega)=\sum_{j=0}^{n-1}(\mu_j-\mu),\ \
\mathcal{H}^*(n,\omega)=\max_{0\le s\le
n-1}\sum_{j=0}^{s}(\mu_j-\mu)
\end{equation}
We shall use the following notations. If $Y$ is a random variable
then $||Y||$ is its usual norm in $L_{2+2\delta}(\Omega)$, that is
\begin{equation}\label{norm}
||Y||=\left(\mathbb{E}|Y|^{2+2\delta}\right)^\frac{1}{2+2\delta},
\end{equation}
where $\delta>0$ is such that $r(2+2\delta)<1$.
\begin{lemma}\label{lem2} In the i.i.d. environment with
$r(2+2\delta)<1$ the following relations hold:
\begin{equation}\label{H*}
 ||\mathcal{H}^*(n)||\le Cn^\frac{1}{2}\ \ \hbox{ and }
\end{equation}
\begin{equation}\label{H2}
\lim_{n\to\infty}n^{-\frac{1+c}{2}} \mathcal{H}(n,\omega)=0 \ \
\hbox{ with $\mathbb{P}$ probability 1 for any $c>0$,}
\end{equation}
The constant $C$ in (\ref{H*}) depends only on $\delta$ and the
distribution of the environment.
\end{lemma}
{\it Remark.} Even though the random variables $\mu_j$ are not
independent, a statement which is stronger than (\ref{H2}) can be
proved. We don't do this because (\ref{H2}) is sufficient for our
purposes.

The proof of Lemma \ref{lem2} will be given at the end of this
section. We now continue the proof of the theorem.

If $\omega$ is such that (\ref{H2}) holds then there is
$t(\omega)$ such that
\begin{equation}\label{Hb}
|b(t,\omega)+\sqrt{t}\sigma^*x - \mu^{-1}t|\equiv
|-\mu^{-1}\mathcal{H}(\mu^{-1}t,\omega)+\sqrt{t}\sigma^*x
|<t^{\frac{1+c}{2}}
 \ \ \hbox{ if $t>t(\omega)$}
\end{equation}
(see (\ref{b})). Hence the following Lemma implies the result we
want:
\begin{lemma}\label{lem3} For a sufficiently small $c>0$
\begin{equation}\label{rel3}
\lim_{n\to\infty}n^{-\frac{1}{2}} \max_{|s|\le
n^{\frac{1+c}{2}}}\left|\sum_{k=n}^{n+s}(\mu_k-\mu)\right|=0 \ \
\hbox{ with $\mathbb{P}$ probability 1}
\end{equation}
\end{lemma}
We put $n=\mu^{-1}t$ in (\ref{Hb});  $\sigma^*$ and $x$ which are
present in (\ref{rel1}) disappear here because of (\ref{Hb}) and the
presence of the small $c$ under the $\max$ sign.

\smallskip\noindent{\it Proof of Lemma \ref{lem3}}. Put
\begin{equation}\label{rel4}
R(n,c,\omega)=n^{-\frac{1}{2}} \max_{|s|\le
n^{\frac{1+c}{2}}}\left|\sum_{k=n}^{n+s}(\mu_k-\mu)\right|
\end{equation}
Note first that if (\ref{rel3}) holds for a subsequences $R(n^2,
\tilde c,\omega)$ then (\ref{rel3}) holds for the whole sequence
$R(n, c,\omega)$, where $c=0.5\tilde c$. Indeed, suppose that
$\omega$ and $\tilde c>0$ are such that (\ref{rel3}) holds for the
subsequence $R(n^2, \tilde c,\omega)$. Then, for $i\in [1, 2n]$, we
have
\begin{equation}\label{rel5}
R(n^2+i,c,\omega)=(n^2+i)^{-\frac{1}{2}} \max_{|s|\le
(n^2+i)^{\frac{1+c}{2}}}\left|\sum_{k=n^2+i}^{n^2+i+s}(\mu_k-\mu)\right|.
\end{equation}
But
\begin{equation}\label{rel6}
\left|\sum_{k=n^2+i}^{n^2+i+s}(\mu_k-\mu)\right|\le
\left|\sum_{k=n^2}^{n^2+i+s}(\mu_k-\mu)\right|+\left|\sum_{k=n^2}^{n^2+i}(\mu_k-\mu)\right|.
\end{equation}
We note next that, since $c<\tilde c$, the inequality $i+|s|<
n^{1+\tilde c}$ holds for sufficiently large values of $n$. This
together with (\ref{rel6}) implies that for sufficiently large $n$
\begin{equation}\label{rel7}
R(n^2+i,c,\omega)\le 2R(n^2,\tilde c,\omega).
\end{equation}
It remains to prove that for that if $\tilde{c}>0$ is small small
enough then
\begin{equation}\label{rel8}
\lim_{n\to\infty}n^{-1} \max_{|s|\le
n^{1+\tilde{c}}}\left|\sum_{k=n^2}^{n^2+s}(\mu_k-\mu)\right|=0 \ \
\hbox{ with $\mathbb{P}$ probability 1}.
\end{equation}
Note that (\ref{H*}) is equivalent to saying that the sequence
$\mu_j$ has the following property: for any $m$
\begin{equation}\label{rel9}
\mathbb{E} \left(\max_{0\le s\le
m}\left|\sum_{k=0}^{s}(\mu_k-\mu)\right|\right)^{2+2\delta}\le
Cm^{1+\delta},
\end{equation}
where $C$ is a constant (related to the previous $C$ in an obvious
way). Using (\ref{rel9}) and stationarity of $\mu_j$ we obtain
that
\begin{equation}\label{rel10}
\mathbb{E}\left(R(n^2,\tilde{c},\omega)\right)^{2+2\delta}\equiv\mathbb{E}
\left(n^{-1} \max_{|s|\le
n^{1+\tilde{c}}}\left|\sum_{k=n^2}^{n^2+s}(\mu_k-\mu)\right|\right)^{2+2\delta}\le
Cn^{(1+\tilde{c})(1+\delta)-2-2\delta}.
\end{equation}
It is now obvious that if  $\tilde{c}<(1+\delta)^{-1}\delta$ then
\begin{equation}\label{rel11}
 \sum_{n=1}^{\infty}\mathbb{E}\left(R(n^2,\tilde{c},\omega)\right)^{2+2\delta}
<\infty
\end{equation}
and the latter in particular implies that
$\lim_{n\to\infty}R(n^2,\tilde{c},\omega)=0$ for almost all
$\omega$. Lemma \ref{lem3} and thus also Theorem \ref{ThMain3} is
proved. $\Box$

\smallskip\noindent{\bf Proof of Theorem \ref{ThMain4}.} In order to check that
(\ref{rel2}) holds we note that $\mu_0$-uniform ergodicity (see
Definition 1) implies that
\begin{equation}\label{rel12}
\left|n^{-1}\sum_{j=k+1}^{k+n}(\mu_0(\mathcal{T}^j\omega)-\mu)\right|\equiv
\left|n^{-1}\sum_{j=k+1}^{n}(\mu_j(\omega)-\mu)\right|\le\varepsilon_n.
\end{equation}
This is due to the fact that, since  $\omega$ in (\ref{uniferg}) is
arbitrary, it can be replaced by $\mathcal{T}^k\omega$. In
particular the left hand side in (\ref{rel2}) can be estimated as
\begin{equation}\label{rel13}
t^{-\frac{1}{2}}
\left|\sum_{k=\mu^{-1}t}^{b(t)+\sqrt{t}\sigma^*x-1}(\mu_k-\mu)\right|\le
\varepsilon_{t^{\frac{1}{2}}}\sigma^*x.
\end{equation}
The proof of (\ref{rel2}) is finished. $\Box$

\smallskip\noindent{\bf Proof of Lemma \ref{lem2}.} It follows from
(\ref{mu}) that
\begin{equation}\label{mu1}
\mu_j-\mu =\sum_{i=1}^\infty B(i,j),\ \ \hbox{ where }\ \
B(i,j)=2(A_j\ldots A_{j-i+1}-r(1)^{i})
\end{equation}
Let us put $\beta=r(2+2\delta)^\frac{1}{2+2\delta}$. Since
$||A_j\ldots A_{j-i+1}||=\beta^i$ and  $r(1)\le \beta $ by
Jensen's inequality, we have
\begin{equation}\label{e1}
||B(i,j)||=\le 2||A_j\ldots A_{j-i+1}||+2r(1)^{i}\le 4\beta^i.
\end{equation}
The $\mathcal{H}(n,\omega)$ can be presented as
\begin{equation}\label{mu2}
\mathcal{H}(n,\omega)=\sum_{j=0}^{n-1}\sum_{i=0}^\infty B(i,j)=
\sum_{i=1}^\infty \sum_{j=0}^{n-1}B(i,j)=\sum_{i=1}^l
\sum_{j=0}^{n-1}B(i,j)+\sum_{i=l+1}^\infty \sum_{j=0}^{n-1}B(i,j),
\end{equation}
where $1\ll l\ll n$ will be chosen later. Denote
\begin{equation}\label{mu3}
B_n(i)= \sum_{j=0}^{n-1}B(i,j)\ \ \hbox{ and }\ \
B_n^*(i)=\max_{0\le s\le n-1}\sum_{j=0}^{s}B(i,j).
\end{equation}
It is then clear that
\begin{equation}\label{H3}
\mathcal{H}^*(n,\omega)\le
\sum_{i=1}^{l}B_n^*(i)+\sum_{i=l+1}^\infty
\sum_{j=0}^{n-1}|B(i,j)|
\end{equation}
and hence
\begin{equation}\label{H4}
||\mathcal{H}^*(n)||\le
\sum_{i=1}^{l}||B_n^*(i)||+\sum_{i=l+1}^\infty
\sum_{j=0}^{n-1}||B(i,j)||\le
\sum_{i=1}^{l}||B_n^*(i)||+4n\frac{\beta^{l+1}}{1-\beta},
\end{equation}
where the last step is due to (\ref{e1}). To estimate
$||B_n^*(i)||$ we note that
\[
B_n(i)=\sum_{j=0}^{n-1}B(i,j)
=\sum_{k=0}^{i-1}\sum_{j=0}^{s_k}B(i,k+ij),\ \hbox{ where }\
s_k=\lfloor\frac{n-k}{i}\rfloor.
\]
Each $D_n(i,k)\overset{\mathrm{def}}{=}\sum_{j=0}^{s_k}B(i,k+ij)$
is a sum of i.i.d. random variables. We put
\[
D_n^*(i,k)\overset{\mathrm{def}}{=}\max_{0\le s\le
s_k}|\sum_{j=0}^{s}B(i,k+ij)|
\]
By Doob's inequality
\[
||D_n^*(i,k)||\le \frac{2+2\delta}{1+2\delta}||D_n(i,k)||
\]
and then by Marcinkiewicz-Zygmund inequality
\[
||D_n^*(i,k)||\le \frac{2+2\delta}{1+ 2\delta}C_\delta
\left(\frac{n}{i}\right)^\frac{1}{2}||B(i,0)||\le
C_1\left(\frac{n}{i}\right)^\frac{1}{2}\beta^i,
\]
where $C_\delta$ depends only on $\delta$ and
$C_1=4\frac{2+2\delta}{1+ 2\delta}C_\delta$. But since
\[
B_n^*(i)\le \sum_{k=0}^{i-1}D_n^*(i,k)
\]
we have that
\[
||B_n^*(i)||\le \sum_{k=0}^{i-1}||D_n^*(i,k)||\le
C_1\left(ni\right)^\frac{1}{2}\beta^i.
\]
Substituting this estimate in (\ref{H4}), we obtain
\begin{equation}\label{H5}
||\mathcal{H}^*(n)||\le
C_1\sum_{i=1}^{l}\left(ni\right)^\frac{1}{2}\beta^i
+4n\frac{\beta^{l+1}}{1-\beta}= C(n,l)n^\frac{1}{2},
\end{equation}
where $C(n,l)=C_1\sum_{i=1}^{l}i^\frac{1}{2}\beta^i
+4n^\frac{1}{2}\frac{\beta^{l+1}}{1-\beta}$. If we now put
$l=n^\frac{1}{2}$, then $\sup_nC(n,n^\frac{1}{2})\le C$ for some
constant $C$. This proves (\ref{H*}). The proof of (\ref{H2})
follows immediately from Lemma \ref{MI} (see Appendix) whose
conditions are satisfied because of (\ref{H*}) and because $\mu_j$
is a stationary sequence. $\Box$
\section{Appendix}
\subsection{Proof of Lemma \ref{simple}.}
Put
\[
r_n(\kappa)\overset{\mathrm{def}}{=}
(\mathbb{E}\prod_{j=1}^nA_j^\kappa)^\frac{1}{n}\ \ \hbox{ and }\ \
f_{n,m}(\kappa)\overset{\mathrm{def}}{=}\max_{0\le s\le m}\log
r_{n+s}(\kappa).
\]
By Jensen's inequality $\log r_n(\kappa)\ge \kappa\mathbb{E}\log
A_1$ and by the same inequality $r_n(\kappa) \le
(\mathbb{E}\prod_{j=1}^nA_j^\gamma)^\frac{1}{n\gamma}=r_n(\gamma)^\frac{1}{\gamma}$.
Condition $\mathbf{C}'$ thus implies that the set of functions
$\{r_n(\cdot)\}$ is uniformly bounded on $[0,\gamma]$ and hence
also the set of functions $\{f_{n,m}(\cdot)\}$ is uniformly
bounded on $[0,\gamma]$. Since functions $r_n(\cdot)$ are convex
on $[0,\gamma]$, the functions $f_{n,m}(\cdot)$ are convex on
$[0,\gamma]$ too. Next,
$f_n(\kappa)\overset{\mathrm{def}}{=}\lim_{m\to\infty}f_{n,m}(\kappa)$
is a limit of functions which converge uniformly on
$[0,\gamma-\epsilon]$, where $\epsilon>0$ is small enough. This
happens because of (a) monotonicity in $m$ of the sequence under
the limit sign, (b) convexity, and (c) existence of bounded right
derivatives $f'_{n,m}(0)$. But then also the monotonically
decaying sequence $f_n(\cdot)$ converges uniformly on
$[0,\gamma-\epsilon]$ (because of the same reasons). Finally,
since $r(\kappa)=\lim_{n\to\infty}f_{n}(\kappa)$, the lemma is
proved. $\Box$

\subsection{Sequences of random variables satisfying the maximal inequality.}
Let $Y_1(\omega), Y_2(\omega),...$ be a sequence of random
variables on the probability space
$(\Omega,\mathcal{F},\mathbb{P})$. We put
\[
S_{k,n}^*(\omega)=\max_{0\le s\le n-1}|\sum_{j=k}^{k+s}Y_j |, \ \
S_{n}(\omega)\mathbb{}=|\sum_{j=1}^{n}Y_j |.
\]
\begin{lemma}\label{MI}
Suppose that for some constant $C$ the inequality $||S_{k,n}^*||\le
Cn^\frac{1}{2}$ holds for all $k,n$. Then
\begin{equation}\label{rel14}
\lim_{n\to\infty}n^{-\frac{1+c}{2}}S_{n}=0 \ \hbox{ with
$\mathbb{P}$ probability 1 for any $c>0$.}
\end{equation}
\end{lemma}
\proof The condition of the lemma implies that
\[
\mathbb{E}(n^{-\frac{1+c}{2}}S_{n})^{2+2\delta}=
n^{-(1+c)(1+\delta)}||S_{n}||^{2+2\delta}\le C_1n^{-c(1+\delta)},
\]
where $C_1=C^{2+2\delta}$. If an integer $m$ is such that
$c(1+\delta)m>1$, then
\[
\sum_{n=1}^{\infty}\mathbb{E}(n^{-\frac{1+c}{2}m}S_{n^m})^{2+2\delta}
\le C_1\sum_{n=1}^{\infty}n^{-c(1+\delta)m}<\infty.
\]
This proves (\ref{rel14}) for the subsequence $n^m$. To control the
rest of the sequence, we shall show that
\[
V(n,l)\equiv |(n^m+l)^{-\frac{1+c}{2}}S_{n^m+l}-
n^{-\frac{1+c}{2}m}S_{n^m}|\to0\ \ \hbox{ as $n\to\infty$ }
\]
uniformly in $l\in[1,(n+1)^m -n^m]$  with $\mathbb{P}$-probability
1. To this end note that
\[
\begin{aligned}
&V(n,l)=|(n^m+l)^{-\frac{1+c}{2}}(S_{n^m+l}-S_{n^m})-
(n^{-\frac{1+c}{2}m}-(n^m+l)^{-\frac{1+c}{2}})S_{n^m}|\\
&\le I_1(n,l)+I_2(n),
\end{aligned}
\]
where $I_1(n,l)=n^{-m\frac{1+c}{2}}|S_{n^m+l}-S_{n^m}|$ and
$I_2(n)=n^{-m\frac{1+c}{2}}S_{n^m}$. We have just proved that
$I_2(n)\to0$ as $n\to\infty$. To estimate $I_1(n,l)$ note that
\[
I_1(n,l)\le n^{-m\frac{1+c}{2}}|\sum_{j=n^m+1}^{n^m+l}Y_j |\le
n^{-m\frac{1+c}{2}}S_{n^m,(n+1)^m-n^m}^*\equiv I_3(n).
\]
But then
\[
\begin{aligned}
& \mathbb{E}(I_3(n))^{2+2\delta}=
n^{-m(1+c)(1+\delta)}||S_{n^m,(n+1)^m-n^m}^*||^{2+2\delta}\\
& \le C_1n^{-m(1+c)(1+\delta)}((n+1)^m-n^m)^{1+\delta}\\
&\le C_2 n^{-m(1+c)(1+\delta)+(m-1)(1+\delta)}= C_2
n^{-(1+\delta)(mc+1)},
\end{aligned}
\]
where the choice of $C_2$ is obvious. It is now clear that
\[
\sum_{n=1}^{\infty} \mathbb{E}(I_3(n))^{2+2\delta}<\infty
\]
and hence $\lim_{n\to\infty}I_3(n)=0$ with $\mathbb{P}$-probability
1. This implies that $I_1(n,l)\to 0$ and thus also $V(n,l)\to 0$ as
$n\to\infty$ uniformly in $l\in[1,(n+1)^m -n^m]$ with
$\mathbb{P}$-probability 1. The lemma is proved. $\Box$

\subsection{General equations for $\mathcal{E}_{x}T_x$ and
$\mathrm{Var}_{x}(T_x)$.}

We shall make use of two general equations. One is the well known
equation for the expectations of hitting times (equation
(\ref{2.1}) below). It can be found in any textbook on Markov
chains. The other (equation (\ref{2.2})) establishes
relationbetween the expectation and the variance of a hitting time
of a random walk. It is equally elementary but it seems that it is
easier to derive it than to find a proper reference. Since the
proof of (\ref{2.2}) naturally includes the derivation of
(\ref{2.1}) both relations are proved here.

Consider a connected Markov chain with a  with a discrete phase
space $S$ and a transition kernel $k(x,y)$, and let $\mathcal{B}$
be a proper subset of $S$. For $x\in S\setminus \mathcal{B}$
denote by $T_x$ the first moment at which the random walk starting
from $x$ hits $\mathcal{B}$. Put
\[
e(x)\overset{\mathrm{def}}{=}\mathcal{E}_{x}(T_x), \ \ \
v(x)\overset{\mathrm{def}}{=}\mathcal{E}_{x}(T_x-e(x))^2\equiv
\mathrm{Var}_{x}(T_x),
\]
where $\mathcal{E}_{x}$ is the usual expectation with respect to
the measure on the space of trajectories starting from $x$. All
expectations considered in this section are supposed to be finite.
\begin{lemma}\label{lemma4} The functions $e(x)$ and $v(x)$
satisfy the following systems of equations:
\begin{equation}\label{2.1}
\left\{
\begin{array}
[c]{ll}%
e(x)=\sum_{y}k(x,y)e(y) +1,& \mathrm{if\quad}x\in S\setminus \mathcal{B},\\
e(x)=0 & \mathrm{if\quad}x\in \mathcal{B},\\%
\end{array}
\right.
\end{equation}
\begin{equation}\label{2.2}
\left\{
\begin{array}
[c]{ll}%
v(x)=\sum_{y}k(x,y)v(y) + f(x),& \mathrm{if\quad}x\in S\setminus \mathcal{B},\\
v(x)=0 & \mathrm{if\quad}x\in \mathcal{B},\\%
\end{array}
\right.
\end{equation}
where $f(x)=\sum_{y}k(x,y)(e(y)-e(x)+1)^2$.
\end{lemma}
\proof Denote by $\chi_{x,y}$ the indicator function of the event
\[
\{\hbox{ the first step of a random walk starting from $x$ is to
$y$}\}.
\]
Obviously $1=\sum_y\chi_{x,y}$ and hence
\begin{equation}\label{2.3}
T_x=\sum_{y}\chi_{x,y}T_x=\sum_{y}\chi_{x,y}(T_y+1).
\end{equation}
Since $\mathcal{E}_{x}\left(\chi_{x,y}(T_y+1)\right)=
k(x,y)(\mathcal{E}_{y}T_y+1)$, applying $\mathcal{E}_{x}$ to both
parts of (\ref{2.3}) leads to (\ref{2.1}).

Similarly
\begin{equation}\label{2.4}
(T_x-e(x))^2=\sum_{y}\chi_{x, y}(T_x-e(x))^2=\sum_{y}\chi_{x,
y}(T_y+1-e(x))^2
\end{equation}
and applying $\mathcal{E}_{x}$ to both parts of (\ref{2.4}) leads
to
\begin{equation}\label{2.5}
v(x)=\sum_{y}k(x,y)\mathcal{E}_{y} (T_y+1-e(x))^2.
\end{equation}
In order to obtain the first equation of (\ref{2.2}) it remains to
observe that
\[
\mathcal{E}_{y} (T_y+1-e(x))^2=\mathcal{E}_{y}
(T_y-e(y))^2+(e(y)-e(x)+1)^2 =v(y)+(e(y)-e(x)+1)^2
\]
and to substitute the last relation into (\ref{2.5}).

Finally, the second equation in (\ref{2.1}) and (\ref{2.2}) is
obvious. $\Box$

\subsection{Proof of Lemma \ref{var}.}

To prove Lemma \ref{var} we shall use the results of the previous
subsection in the case when $S=\mathbb{Z}$ is a line and
$\mathcal{B}\equiv \mathcal{B}_n$ is a semi-line of integers which
are $\ge n$. Technically, equations (\ref{2.1}) are a particular
case of (\ref{2.2}) and it makes sense to solve that latter for a
general function $f(x)$. Note first that (\ref{2.2}) can be
re-written in terms of parameters $p_{k},\ -\infty<k<\infty,$ as
follows:
\begin{equation}
\left\{
\begin{array}
[c]{ll}%
g_k=p_{k}g_{k+1}+ q_{k}g_{k-1}+f_k& \mathrm{if\quad}k<n,\\
g_n=0 ,& %
\end{array}
\right. \label{4.1}%
\end{equation}
where the meaning of $g_k$ depends on the choice of $f$. Solving
(\ref{4.1}) is a relatively simple and well studied matter. The
following lemma is included into this work for the sake of
completeness. As before, $A_j=p_jq_j^{-1}\equiv p_j(1-p_j)^{-1}$;
the sequence $\omega=(p_j)_{-\infty<j<\infty}$ is fixed throughout
this section.

\begin{lemma}\label{lemma5} Suppose that

(i) $\sum_{j=0}^{\infty} \prod_{i=0}^jA_{-i} <\infty$ and

(ii) $f_k$ is such that
$\sum_{j=0}^{\infty}|f_{-j}|\prod_{i=0}^jA_{-i} <\infty$.

\noindent then the solution  $(g_k),\ -\infty<k\le n-1$, to
(\ref{4.1}) is given by
\begin{equation}\label{4.7}
g_k=\sum_{j=k}^{n-1}d_j,\ \hbox{ where }\ d_j=\sum_{i=0}^\infty
A_j...A_{j-i+1}p_{j-i}^{-1}f_{j-i}.
\end{equation}
This solution can be obtained as $g_k=\lim_{a\to-\infty}h_k $,
where $h_k$ is a solutions to
\begin{equation}
\left\{
\begin{array}
[c]{ll}%
h_k=p_{k}h_{k+1}+ q_{k}h_{k-1}+f_k& \mathrm{if\quad}a<k<n,\\
h_a=h_n=0 ,& %
\end{array}
\right. \label{4.2}%
\end{equation}
\end{lemma}

\proof To solve (\ref{4.2}), present $h_k$ as
\begin{equation}\label{4.1.1}
h_{k}=\varphi_{k}h_{k+1}+\tilde{d}_{k},\ \ k\ge a.
\end{equation}
If we put $\varphi_a=0$ and $\tilde{d}_a=0$, then an easy
induction argument (involving (\ref{4.2})) leads to the following
formulae:
\begin{equation}\label{4.11}
 \varphi_{k}=(1-q_k\varphi_{k-1})^{-1}p_k,\ \ k\ge a+1
\end{equation}
\begin{equation}\label{4.1.2}
\tilde{d}_k=A_k\tilde{d}_{k-1}+w_k,\ \ k\ge a+1,\ \ \hbox{ where
$w_k=(1-q_k\varphi_{k-1})^{-1}f_k$}.
\end{equation}
Iterating (\ref{4.1.1}) and (\ref{4.1.2}) leads to
\[
h_k=\tilde{d}_k+\varphi_{k}\tilde{d}_{k+1}+...
+\varphi_{k}\varphi_{k+1}...\varphi_{n-1}\tilde{d}_{n-1}
\]
and
\[
\tilde{d}_k=w_k+A_kw_{k-1}+...+A_k...A_{a+2}w_{a+1}.
\]
It follows from (\ref{4.11}) that $0\le \varphi_k<1$ and (direct
calculation)
$1-\varphi_k=q_k(1-q_k\varphi_{k-1})^{-1}(1-\varphi_{k-1})$. Hence
\[
1-\varphi_k\le A_k(1-\varphi_{k-1})\le A_kA_{k-1}\ldots
A_{a+1}\to0 \hbox{ as }\ a \to-\infty,
\]
where the last relation follows from condition $(i)$ of the Lemma.
In other word, $\lim_{a\to-\infty} \varphi_k=1$, and condition
$(ii)$ now implies that $\lim_{a \to-\infty}\tilde{d}_k=d_k$ and
hence $\lim_{a \to-\infty}\tilde{h}_k=g_k$. $\Box$

We shall now prove Lemma \ref{var}. To this end note first that if
we substitute $f_k\equiv1$ into (\ref{4.7}) and (\ref{4.2}), then,
according to (\ref{2.1}) we obtain formulae for $e_k\equiv
\mathcal{E}_\omega T_k(n)$ and thus also for $\mu_k=e_{k+1}-e_{k}$
(see Lemma \ref{Etau}). Next, to find $v_k\equiv\Vo T_k(n)$ we
have to put
\[
f_k=p_k(e_{k+1}-e_{k}+1)^2+q_k(e_{k-1}-e_{k}+1)^2=
p_k(\mu_k+1)^2+q_k(1-\mu_{k-1})^2.
\]
The main equation in (\ref{4.7}) can be rewritten as
\[
p_k(g_k-g_{k+1})=q_k(g_{k-1}-g_k)+f_k\ \ \hbox{ and thus }\ \
(g_k-g_{k+1})=A_k(g_{k-1}-g_k)+p_k^{-1}f_k.
\]
In particular this leads to the following relations:
\[
\mu_k=A_k\mu_{k-1}+p_k^{-1}.
\]
To see now that $d_j$ in (\ref{4.7}) turns into (\ref{2.12}) a
matter of very simple calculation.

Relation (\ref{2.13}) follows now from the condition $r(2)<1$.
Finally the explicit expression (\ref{2.14}) is again a matter of
simple calculation. Lemma \ref{var} is proved. $\Box$

\end{document}